\documentclass{article}[12pt]
\usepackage{amsfonts}
\usepackage{amsmath}
\usepackage{amssymb}
\usepackage{amscd}
\usepackage{color}
\usepackage{amsthm}
\usepackage{indentfirst}
\usepackage[hmargin=3cm,vmargin=3cm]{geometry}

\newtheorem{thm}{Theorem}

\newtheorem{lma}{Lemma}

\theoremstyle{definition}

\theoremstyle{remark}

\newcommand{\R}{{\mathbb{R}}}

\newcommand{\C}{{\mathbb{C}}}

\newcommand{\D}{{\mathbb{D}}}

\newcommand{\del}{\partial}

\newcommand{\zbar}{\overline{z}}

\newcommand{\G}{\mathcal{G}}

\newcommand{\om}{\omega}

\begin{document}

\title{\textbf{Enlacements asymptotiques revisit\'es.}}

\author{\textsc{Egor Shelukhin.}
\\
}

\date{}
\maketitle

\begin{abstract}
We give an alternative proof of a theorem of Gambaudo-Ghys \cite{Enlacements} and Fathi \cite{FathiThese} on the interpretation of the Calabi homomorphism for the standard symplectic disc as an average rotation number. This proof uses only basic complex analysis.
\end{abstract}

\tableofcontents

\section{The theorem of Gambaudo-Ghys and Fathi.}

Let $\G = Ham_c(\D,\om)$ be the group of compactly supported Hamiltonian diffeomorphisms of the standard disc $\D = \{z \in \C \; : \; |z|\leq 1\}$ endowed with the standard symplectic form $\om = \frac{i}{2}dz \wedge d\overline{z}.$
The Calabi homomorphism \cite{CalabiHomomorphism} from $\G$ to $\R$ is defined as \[Cal(\phi) = \int_0^1 dt \int_{\D} H_t \om,\] where $H_t$ is the normalized Hamiltonian (zero near the boundary) of a Hamiltonian isotopy $\{\phi_t\}_{t\in [0,1]}$ with $\phi_1 = \phi.$ In other words this isotopy is generated by a time-dependent vector field $X_t,$ that satisfies the relation \[\iota_{X_t} \om = - d H_t.\]

The mean rotation number is defined in terms convenient for the proof as follows. Consider the differential form \[\alpha =\frac{1}{2\pi}\frac{d(z_1 -z_2)}{z_1 - z_2}\] (used by Arnol'd \cite{ArnoldBraids} in his study of the cohomology of the pure braid groups) on the configuration space $X_2 = X_2(\D) = \{(z_1,z_2)| z_j \in \D, \; z_1 \neq z_2\} = \D \times \D \setminus \Delta,$ where $\Delta \subset \D \times \D$ is the diagonal. Denote by \[\theta = Im(\alpha)\] its imaginary part. Note that the two forms $\alpha$ and $\theta$ are closed. For each pair of points $(z_1,z_2) \in \D \times \D$ such that $z_1 \neq z_2,$ that is for each point $x = (z_1,z_2) \in X_2,$ consider the curve $\{\phi_t \cdot x\}$ in $X_2$ defined by \[\phi_t \cdot x = (\phi_t(z_1), \phi_t(z_2))\] for each $t \in [0,1].$ The average rotation number is \[\Phi(\phi) = \int_{X_2} dm^2(x) \int_{\{\phi_t \cdot x\}} \theta,\] where $dm^2(x) = dm(z_1)dm(z_2)$ is the Lebesgue measure on $\D \times \D$ restricted to $X_2.$ By preservation of volume, it is clear that $\Phi$ is a homomorphism $\G \to \R.$

The theorem of Gambaudo-Ghys \cite{Enlacements} and Fathi \cite{FathiThese} is the following equality.

\begin{thm}\label{GG theorem}
\[\Phi = -2Cal,\]
as homomorphisms $\G \to \R.$

\end{thm}

Gambaudo and Ghys have presented several proofs of this result, and in \cite{FathiThese} a different proof of Fathi is found. More proofs of this result are known today (cf. \cite{Deryabin}). Here we present an alternative short proof, which is in fact a complex-variable version of the proof of Fathi.

\subsection*{Acknowledgements}
I thank Steven Lu for inviting me to give a talk on the CIRGET seminar in Montr\'eal, that has lead me to revisit the theorem Gambaudo-Ghys and Fathi. I thank Albert Fathi for sending me his original proof of the theorem. I thank Boris Khesin for the reference \cite{Deryabin}.

\section{The alternative proof.}

Put $\xi_t = dz(X_t),$ for the natural complex coordinate $z$ on $\D.$ Hence $\xi_t$ is a smooth complex-valued function on $\D,$ vanishing near the boundary. The computations \[\iota_{X_t} (\frac{i}{2} dz\wedge d\overline{z}) = \frac{i}{2} \xi_t d\overline{z} - \frac{i}{2} \overline{\xi}_t dz\] and \[-dH_t = - \frac{\partial H_t}{\partial \overline{z}} d\zbar - \frac{\partial H_t}{\partial z} dz\] give us \begin{equation}\label{xi and del H del zbar}\xi_t = 2i \frac{\del H_t}{ \del \zbar}.\end{equation}

Now \[\Phi(\phi) = Im(\int_{X_2} dm^2(x) \int_{\{\phi_t \cdot x\}} \alpha),\] and hence it is sufficient to compute

\[\int_{X_2} dm^2(x) \int_{\{\phi_t \cdot x\}} \alpha = \frac{1}{2\pi} \int_{X_2} dm^2(x) \int_{\{\phi_t \cdot x\}} \frac{d(z_1 - z_2)}{z_1 - z_2} =\]

\[= \frac{1}{2\pi} \int_{X_2} dm(z_1)dm(z_2)\int_0^1 dt \; \frac{\xi_t(\phi_t(z_1)) - \xi_t(\phi_t(z_2))}{\phi_t(z_1) - \phi_t(z_2)} = \]

as the function is absolutely integrable (see Lemma \ref{integrabilite absolument}), by Fubini,

\[ =\frac{1}{2\pi} \int_0^1 dt \; \int_{X_2} dm(z_1)dm(z_2) \; \frac{\xi_t(\phi_t(z_1)) - \xi_t(\phi_t(z_2))}{\phi_t(z_1) - \phi_t(z_2)} = \frac{1}{2\pi} \int_0^1 dt \; \int_{X_2} dm(z_1)dm(z_2) \; \frac{\xi_t(z_1) - \xi_t(z_2)}{z_1 - z_2} = \]

as both terms of the sum are absolutely integrable (a consequence of the proof of Lemma \ref{integrabilite absolument} as well),

\[ = 2 \cdot \frac{1}{2\pi} \int_0^1 dt \; \int_{\D} dm(w) \int_{\D \setminus \{w\}} \frac{\xi_t(z)}{z - w} dm(z) = \]

by Equation \ref{xi and del H del zbar},

\[ = - 2 i \int_0^1 dt \; \int_{\D} dm(w) \int_{\D \setminus \{w\}} \frac{1}{2 \pi i}\frac{\del H_t}{\del \zbar}\frac{dz \wedge d\zbar }{z - w} = - 2 i \int_0^1 dt \; \int_{\D} dm(w) H_t(w) = - 2 i Cal(\phi).\]

The penultimate equality is a consequence of the Cauchy formula for smooth functions \cite[Theorem 1.2.1]{HormanderComplex}. For any $C^1$ function $f:\D \to \C,$ we have \[f(w) = \frac{1}{2 \pi i} \int_{\partial \D} \frac{f(z)}{z-w} dz + \frac{1}{2\pi i} \int_{\D} \frac{\del f}{\del \zbar} \frac{dz \wedge d\zbar}{z-w}.\] It remains to note that as $H_t$ is zero near the boundary, the first term of the sum vanishes.

Now we show the absolute integrability that we use to change the order of integration.

\begin{lma}\label{integrabilite absolument}

\[\int_{X_2}\int_0^1 dm(z_1)dm(z_2) dt \; \frac{|\xi_t(\phi_t(z_1)) - \xi_t(\phi_t(z_2))|}{|\phi_t(z_1) - \phi_t(z_2)|} < \infty \]

\end{lma}

By the Tonnelli theorem, the following chain of inequalities suffices:

\[\int_0^1 dt \int_{X_2} dm(z_1)dm(z_2) \; \frac{|\xi_t(\phi_t(z_1)) - \xi_t(\phi_t(z_2))|}{|\phi_t(z_1) - \phi_t(z_2)|} = \]

\[= \int_0^1 dt \int_{X_2} dm(z_1)dm(z_2) \; \frac{|\xi_t(z_1) - \xi_t(z_2)|}{|z_1 - z_2|} \leq 2 \int_0^1 dt \int_{z \in \D} dm(z) |\xi_t(z)| \int_{w \in \D \setminus \{z\}}  \; \frac{1}{|z - w|} dm(w) \leq \]

\[ \leq 8\pi \int_0^1 dt \int_{\D} |\xi_t| dm < \infty,\]

because \[\int_{w \in \D \setminus \{z\}}  \; \frac{1}{|z - w|} dm(w) \leq 4\pi,\] as one verifies by direct calculation.

\bibliographystyle{amsplain}

\vspace{3mm}

\textsc{Egor Shelukhin, D\'{e}partement de math\'{e}matiques et de Statistique, Universit\'{e} de Montr\'{e}al, C.P. 6128, Succ. Centre-
ville, Montr\'{e}al H3C 3J7, Qu\'{e}bec, Canada}\\

\end{document}